\documentstyle[12pt]{article}

  \newcommand{\const}{\rm const}
  \newcommand{\Var}{\rm Var}
  
  \newcommand{\Law}{\rm Law}

  \newcommand{\Conv}{\rm Conv}

\begin{document}

   \begin{center}

 {\bf  Generalization  and refinement  of Khintchin's inequality.}\par

\vspace{5mm}

{\bf M.R.Formica,  E.Ostrovsky and L.Sirota}.\par

\end{center}

\vspace{5mm}

 Universit\`{a} degli Studi di Napoli Parthenope, via Generale Parisi 13, Palazzo Pacanowsky, 80132,
Napoli, Italy. \\

e-mail: mara.formica@uniparthenope.it \\

\vspace{4mm}

Department of Mathematics and Statistics, Bar-Ilan University, \\
59200, Ramat Gan, Israel. \\

e-mail: eugostrovsky@list.ru\\

\vspace{4mm}

\ Department of Mathematics and Statistics, Bar-Ilan University, \\
59200, Ramat Gan, Israel. \\

e-mail: sirota3@bezeqint.net \\

\vspace{5mm}

\begin{center}

{\bf Abstract.} \\

\end{center}

\vspace{4mm}

 \ We derive the exponential as well as power decreasing  tail estimations for normed sums of centered
 independent identical distributed (or not) random variables on the Khintchine's  form.\par
 \ We consider arbitrary, in particular, non - Rademacher's  variables and not only Lebesgue - Riesz rearrangement invariant  norms for
 the random variables.\par
 \ We intend to calculate the {\it exact } value of correspondent limit. \par

\vspace{4mm}

\ {\it Key words and phrases.} Probability,  random variables (r.v.) and sequences, independent (i) and identical distributed (i.d.) r.v.,
 tail of distribution, moment,  normalized sums, generalized  and ordinary Khintchin's form and constants, convex functions, Cramer's condition,
 Young - Fenchel, or Legendre transform, generating and moment generating function,  subgaussian variable and condition,
 Lebesgue - Riesz, Orlicz and Grand Lebesgue Spaces (GLS) and norms,
 Rademacher's r.v.,  examples, martingales, subgaussian variables, Pythagoras inequality, Khintchine's statement of problem and inequality. \par

\vspace{5mm}

\section{Statement of problem. Notations. Previous results.}

\vspace{5mm}

 \hspace{3mm} Let $ \ ( \ \Omega = \{\ \omega \ \}, \cal{B}, {\bf P} \ )  \ $ be certain probability space with expectation $ \ {\bf E} \ $
 and variance $ \ \Var. \ $ It will be presumed that it is sufficiently rich; so that there exists considered further random sequences. \par

  \ The ordinary Lebesgue - Riesz $ \ L_p = L_p(\Omega) \ $ for the numerical valued random variable $ \ \eta \ $
 is defined as usually

$$
||\eta||_p \stackrel{def}{=} \left[ \ {\bf E} |\eta|^p \ \right]^{1/p}, \ p \ge 1.
$$

  \ Let also $ \ \{\xi_i \}, i = 1,2, \ldots \ $ be a sequence of centered (mean zero) independent identical distributed (i., i.d.)
 random variables (r.v.); $ \ \xi = \xi_1. \ $ The law of distribution of the r.v. $ \ \xi = \xi_1 \ $ will be denoted by $ \ L(\xi). \ $ \par

 \  Denote by $ \ D(n), \ n = 1,2,\ldots \ $ the $ \ n \ - \ $  tuple  (vector) of the deterministic  numerical sequences of the form

$$
D(n) \stackrel{def}{=} \{ \ (a_1,a_2,  \ldots, a_n)  \ \} = \{ \ \vec{a} \ \},
$$
for which

$$
[\sum_{k=1}^n  a^2_k]^{1/2} = ||\vec{a}||_2 = ||a||_2 = 1,
$$
and put

$$
 ||a||_p = ||\vec{a}||_p   \stackrel{def}{=} [\sum_{k=1}^n |a_k|^p]^{1/p}, \ p \ge 2;
$$

$$
D := \cup_{n=1}^{\infty} D(n);
$$

$$
S_n = S_n[\xi] = S_n[D(n)] = S_n[D(n)](\{\xi_k\}) \stackrel{def}{=} \frac{\sum_{k=1}^n a_k \xi_k}{\sqrt{\sum_{k=1}^n a^2_k}} =\sum_{k=1}^n a_k \xi_k,
$$
as long as {\it we assume henceforth}

\begin{equation} \label{cond s Dn}
a = \vec{a} = \ (a_1,a_2,  \ldots,a_n) \in D(n).
\end{equation}

\vspace{4mm}

 \ {\bf Definition 1.1.}  \  Recall the definition of the following important {\it Khintchin's constants}

\begin{equation} \label{Bp constants}
B[L(\xi)](p) \stackrel{def}{=} \sup_{n = 1,2,\ldots} \ \sup_{ \{a_k\} \in D(n)} \frac{ \ ||S_n[D(n)](\{\xi_k\}) \ ||_p}{\sqrt{\sum_{k=1}^n a^2_k}},
\end{equation}

\vspace{3mm}

\begin{equation} \label{Bp constants}
A[L(\xi)](p) \stackrel{def}{=} \inf_{n = 1,2,\ldots} \ \inf_{ \{a_k\} \in D(n)} \frac{ \ ||S_n[D(n)](\{\xi_k\}) \ ||_p}{\sqrt{\sum_{k=1}^n a^2_k}}.
\end{equation}

  \ These constants was introduced by A.Ya.Khintchin in \cite{Khintchin} for the Rademacher's distribution of the (independent) r.v. - s $ \ \{\xi_k\} \ $ and
for the Lebesgue - Riesz spaces $ \ L_p = L_p(\Omega, {\bf P}): \ $

$$
{\bf P}(\xi_k = 1) = {\bf P}(\xi_k = -1) = 1/2.
$$
 \hspace{3mm} They was investigated in many works: \cite{Dilworth}, \cite{Figel}, \cite{Haagerup},  \cite{Kahane1}, \cite{Kahane2},
 \cite{Khare}, \cite{Komorowski}, \cite{Kwapien}, \cite{Latala}, \cite{Montgomery-Smith}, \cite{Peskir}, \cite{Rodin}, \cite{Serb},
 \cite{Szarek}, \cite{Tomaszewski1}, \cite{Tomaszewski2} etc. In particular, was obtained the exact values of these constants as well as its
 upper and lover bounds, as a rule, for the Rademacher's variables. \par
  \ It was considered also the case when the r.v. $ \ \xi, \ $ as well as its copies,  belongs to some separable Banach space.\par

\vspace{5mm}

\hspace{3mm} {\bf  The aim of this preprint is twofold: a generalization of the mentioned results into the case of an arbitrary distribution of the
r.v. - s  and into the another rearrangement invariant Banach spaces, instead the classical Lebesgue - Riesz ones, in particular, into the
 Grand Lebesgue Spaces (GLS), builded on the source probability space.} \par

{\bf \ We clarify also  the known results and obtain sometimes the exact values of correspondent constants, see for example
 theorem 3.1.} \par

\vspace{4mm}

 \ The case when the r.v. - s $ \ \{\xi_k\} \ $ forms a sequence of martingale differences was investigated in \cite{Osekovski},
\cite{Ostrovsky 2004},  \cite{Ostrovsky 2012}, \cite{Ostrovsky 2014}, \cite{Peskir},  \cite{Peskir2}. Note that the correspondent
Khintchine estimate is formulated in  slightly different terms.  \par

\vspace{5mm}

 \ In detail, let $ \ F \ $ be some Banach rearrangement invariant functional space  builded on the source probability space equipped with the norm  $ \ ||\eta||F \ $
 on the  centered r.v. $ \ \eta, \ $  and let $ \ \eta_j, \ j = 1,2,\ldots \ $ be independent copies of $ \ \eta. \ $ \par

 \vspace{4mm}

\ {\bf Definition 1.2.}  \ Define the following important  for us {\it generalized Khintchin's constants}

\begin{equation} \label{F upper constants}
B[L(\eta)]\{F\} \stackrel{def}{=} \sup_{n = 1,2,\ldots} \ \sup_{ \{a_k\} \in D(n)} \frac{ \ ||S_n[D(n)](\{\eta_k\}) \ ||F}{\sqrt{\sum_{k=1}^n a^2_k}},
\end{equation}

\vspace{3mm}

\begin{equation} \label{F lower constants}
A[L(\eta)]\{F\} \stackrel{def}{=} \inf_{n = 1,2,\ldots} \ \inf_{ \{a_k\} \in D(n)} \frac{ \ ||S_n[D(n)](\{\eta_k\}) \ ||F}{\sqrt{\sum_{k=1}^n a^2_k}}.
\end{equation}

  \ These constants may be named as {\it  generalized } Khintchine's constants defined for the r.v. - s $ \ \{\eta_j \} \ $  and for the  space $ \ F. \ $ \par

\vspace{4mm}

 \hspace{3mm} {\bf Proposition 1.1. Preliminary bounds.} \
 \  Note for the beginning that always $ \ B[L(\eta)]\{F\} \ge ||\eta||F. \ $ Further, denote $ \  \sigma^2 := \Var (\eta) \ $ and
 introduce also the Gaussian distributed centered r.v. $ \ \zeta \ $ with parameters $ \ (0,\sigma^2): \ \Law (\zeta) = N(0,\sigma^2), \ $
 if of course $ \ \sigma^2 \in (0,\infty). \ $ \par

 \ It follows from the classical CLT that $ \  B[L(\eta)]\{F\} \ge ||\zeta||F. \ $ Thus,

 \begin{equation} \label{upp B}
  B[L(\eta)]\{F\} \ge \max(||\zeta||F, \ ||\eta||F).
 \end{equation}

 \ Quite analogously

 \begin{equation}  \label{lower A}
  A[L(\eta)]\{F\} \le \min(||\zeta||F, \ ||\eta||F).
 \end{equation}

 \vspace{3mm}

  \ Notice that the best values of the Khinchine's constant for the classical  Rademacher's series (on the real line)
  of the r.v. - s $ \ \{ \ \xi_j \} \ $ was found by U.Haagerup in \cite{Haagerup}; in this case both the estimations (\ref{upp B})  and (\ref{lower A})
 are exact for the Lebesgue - Riesz spaces $ \ F = L_p \ $  for all the  greater values $ \ p. \ $ \par

  \vspace{5mm}

  \section{Brief description of the theory of Grand Lebesgue Spaces (GLS), with addition. } \par

 \vspace{5mm}

 \hspace{3mm} We will deal with the so - called  Grand Lebesgue Spaces (GLS).  Recall briefly some  used further facts from the theory
 of these spaces. \par

 \vspace{4mm}

  \  {\bf Definition   2.1,}  see   \cite{Kozachenko1},   \cite{Ostrovsky1}, chapter 1, sections 1.1 - 1.3;
 \cite{Ahmed Fiorenza Formica at all}, \cite{fiokarazanalanwen2004}, \cite{fioguptajainstudiamath2008},
\cite{Fiorenza-Formica-Gogatishvili-DEA2018},
\cite{fioforgogakoparakoNA}, \cite{fioformicarakodie2017},
\cite{formicagiovamjom2015}. Let  $ \ \phi = \phi(\lambda). \
\lambda \in (-\lambda_0,\lambda_0),  \lambda_0 = \const \in (0,
\infty] \ $ be even twice  continuous differentiable convex
function, such that $ \ \lambda \to 0 \ \Rightarrow \phi(\lambda)
\asymp \lambda^2, \ $ strictly increasing on the interval $ \ [0,
\lambda_0). \ $   We impose the  following condition on these
functions in the case when $ \ \lambda_0 = \infty: \ $

$$
\lim_{\lambda \to \infty} \frac{\phi(\lambda)}{\lambda} = \infty.
$$

 \ The set of all such a functions will be denoted by $ \ \Phi, \  \Phi = \{ \ \phi \ \}. \ $ \par

 \ By definition,  the random variable  $  \ \zeta \ $ belongs to the space  $  \ B(\phi), \ $  for certain
fixed function $  \phi \in \Phi, $ if and only if there exists a non - negative {\it finite} constant  $ \ \tau \ $ such that

\begin{equation} \label{Bphi}
\forall  \lambda: \ |\lambda| < \lambda_0 \ \Rightarrow {\bf E} \exp(\lambda \zeta) \le \exp(\phi(\lambda \tau)).
\end{equation}

 \ The minimal value of the constant  $  \ \tau \ $ which satisfies the inequality   (\ref{Bphi}) is said to be the
$  B(\phi) \  $ norm  of the r.v. $  \ \zeta: $

$$
||\zeta||B(\phi) \stackrel{def}{=} \max_{\pm} \sup_{\lambda \in (0, \lambda_0)} \phi^{-1} \{ \ln {\bf E} \exp ( \pm \lambda \ \zeta) \}/|\lambda|,
$$
so that

\begin{equation} \label{key ineq}
\forall \lambda: |\lambda| < \lambda_0 \ \Rightarrow  {\bf E} \exp(\lambda \zeta) \le \exp ( \phi(\lambda ||\zeta||B(\phi) ) ).
\end{equation}

  \ We suppose in fact that the  r.v. $ \ \zeta \ $ satisfies the well - known Cramer's
condition:

$$
\exists  c > 0 \ \Rightarrow {\bf P} (|\zeta| > x) \le \exp (-cx), \ x \ge 0.
$$
 then $ \ \lambda_0 > 0. \ $ In this case the generated function $  \ \phi(\cdot) $   may be introduced naturally. Namely, the
 so - called natural  (generating) function  for the r.v. $ \ \zeta \ $ satisfying the Cramer's condition is defined as follows

$$
\phi_{\zeta}(\lambda) := \max_{\pm} \ln {\bf E} \exp(\pm \lambda \ \zeta).
$$
 \ The function $ \ \lambda \to {\bf E} \exp (\lambda \zeta), \ \lambda = \const   \ $  is said to be a {\it moment generating function}
for the r.v. $ \ \zeta, \ $ of course, if there exists in certain non - zero neighborhood of origin. \par
 \ The natural function  $ \ \phi_{\zeta}(\lambda) \ $ play a very important role, in particular, in the theory of Large Deviations (L.D.)

\vspace{4mm}

 \ These $ \ B(\phi) \ $ spaces are complete Banach functional and rearrangement invariant, as
well as considered further Grand Lebesgue Spaces. They were introduced at first
in the article   \cite{Kozachenko1}.  The detail investigation of these spaces may be found in the
the monographs  \cite{Buldygin} and   \cite{Ostrovsky1}, chapters 1,2. \par

 \ It is known that  $  \ 0 \ne \zeta \in B(\phi) $ if and only if $  \ {\bf E} \zeta = 0 $ and

$$
\exists K = \const \in (0, \infty) \ \Rightarrow \max \left[  {\bf P} (\zeta \ge u), {\bf P}(\zeta \le - u) \right] \le
\exp \left\{ - \phi^*(u/K)   \right\}, \ u > 0,
$$
where $ \ \phi^*(u) \ $ denotes the famous Young - Fenchel, or Legendre transform of the function $ \ \phi: \ $

$$
\phi^*(u) \stackrel{def}{=} \sup_{|\lambda| < \lambda_0} (\lambda u - \phi(\lambda));
$$

and herewith

$$
||\zeta||B(\phi) \le C_1(\phi) K \le C_2(\phi) ||\zeta||B(\phi).
$$

 \ More exactly, if $  0 < ||\zeta||B(\phi) = ||\zeta|| < \infty, $ then $ \ \forall u \ge 0 \ \Rightarrow \ $

$$
\max \left[ {\bf P}(\zeta \ge u), \ {\bf P}(\zeta \le - u)   \right] \le \exp \left(  - \phi^*(u/ ||\zeta||)  \right).
$$

\vspace{4mm}

 \ Define the following Young - Orlicz $ \ N - \ $ function

\begin{equation} \label{Orlicz}
N[\phi](u):= \exp\phi^*(u) - 1.
\end{equation}

 \ It is proved in particular in  \cite{Kozachenko1}, \cite{Ostrovsky1}, chapter 1, that if $ \ \lambda_0 = \infty, \ $
then the space $ \ B(\phi) \ $ coincides up to norm equivalence   with the closed subspace of  {\it exponential} Orlicz space $ \ L(N[\phi]) \ $
builded on the source probability space and consisted only on the centered random variables. \par

\vspace{4mm}

 \ Recall yet, see e.g.  \cite{Ahmed Fiorenza Formica at all}, \cite{fiokarazanalanwen2004}, \cite{fioguptajainstudiamath2008},
\cite{Fiorenza-Formica-Gogatishvili-DEA2018},
\cite{fioforgogakoparakoNA}, \cite{fioformicarakodie2017},
  that the so - called Grand Lebesgue Space
  (GLS)  $ \ G\psi \ $  equipped with the norm  $ \ ||\zeta|| G\psi \ $  of the r.v. $ \ \zeta \ $ is defined as follows

$$
||\zeta||G\psi \stackrel{def}{=} \sup_{p \ge 2}  \left[ \ \frac{||\zeta||_p}{\psi(p)} \ \right].
$$
 \ Here $ \ \psi = \psi(p), \ p \ge 1 \ $ is measurable bounded from below  function, which is names as ordinary as
 {\it generating function} for this space. \par
 \ If the r.v. $  \zeta $  belongs to some $ \ B(\phi) \ $ space, then it belongs also to certain  $ \ G\psi \ $ space with

$$
\psi = \psi_{\phi}(p) = \frac{\phi^{-1}(p)}{p}, \ p \ge 2.
$$

 \ The inverse conclusion in not true. Namely, the mean zero r.v. $ \zeta $  can has finite
all the moments $ \ |\zeta|_p < \infty, \ p \ge 2,  \ $ but may not satisfy the Cramer's condition. \par
 \ A very popular class of these spaces form the  {\it subgaussian} random variables, i.e. for
which  $  \ \phi(\lambda) = \phi_2(\lambda) \stackrel{def}{=} 0.5 \lambda^2 $ and $  \lambda_0 = \infty. \ $  For instance, every
centered Gaussian distributed r.v. is also subgaussian. \par
 \  The correspondent   $  \ \psi \ $ function has a form  $ \ \psi(p) = \psi_2(p) = \sqrt{p}. $ \par

 \ More generally, suppose

$$
 \phi(\lambda) = \phi_m(\lambda) = |\lambda|^m/m, \  |\lambda| \ge 1, \ \lambda_0 = \infty, \ m = \const \ge 1.
$$
 \ The correspondent    $ \psi $  function has a form

$$
\psi(p) = \psi_m(p) = p^{1/m}
$$
 and the correspondent tail estimate is follow:

$$
\max \left[  {\bf P}(\zeta \ge u), \ {\bf P} (\zeta \le - u)  \right] \le \exp \left\{ -  (u/K)^m   \right\}, \ u > 0.
$$

 \ These space are used   in particular for obtaining of the exponential estimates for sums of independent random variables, see e.g.
\cite{Kozachenko1},  \cite{Ostrovsky1},  sections 1.6, 2.1 - 2.5. Indeed, introduce for any
function $ \ \phi(\cdot) \ $ from the set $ \ \Phi \ $
a new function   $   \  \overline{\phi}(\cdot)   \  $  which belongs also at the same set:

$$
 \overline{\phi}(\lambda)  \stackrel{def}{=} \sup_{n = 1,2, \ldots} n \  [ \ \phi \{\lambda/\sqrt{n} \} \ ].
$$
 \ It is easily to see that

$$
\sup_n {\bf E} \exp(\lambda S(n)/\sqrt{n}) \le \exp [ \overline{\phi}(\lambda)]
$$
with correspondent uniform relative the variable $  \  n  \ $ exponential tail  estimate. \par

 \ For instance, if $ \ {\bf E} \xi = 0 \ $  and  for some value $  \ m = \const > 0 $

$$
\max \left[  {\bf P}(\xi \ge u), \ {\bf P} (\xi \le - u)  \right] \le \exp \left\{ -  u^m   \right\}, \ u > 0,
$$
then

$$
 \sup_n \max \left[  {\bf P}(S(n)/\sqrt{n} \ge u), \ {\bf P} (S(n)/\sqrt{n} \le - u)  \right] \le
$$

$$
\exp \left\{ - C(m)  u^{\min(m,2)} \right\}, \ u > 0, \ C(m) \in (0, \infty),
$$
and the last estimate is essentially non - improvable. \par

\vspace{4mm}

 \ We will use the following fact, see \cite{Ostrovsky1}, chapter 1, section 1.6, theorem 1.6.1. Let $ \ \{\xi_i\}, \ i = 1,2,\ldots,n \ $ be
 independent r.v.  belonging to the certain space $ \ B(\phi). \ $   Denote the set of all convex functions by $ \ \Conv. \ $ \par
 \  Let us impose the following  important condition on this function ({\sc Condition Triangle })  of order  $ \ r,  r = \const \in [1,2]. \ $ \par

 \vspace{4mm}

 \ {\bf Definition 1.4.,} see \cite{Ostrovsky1} chapter 1, section 1.6.  We will write that  the function $ \ \phi(\cdot) \ $ from the set $ \ \Phi \ $
 belongs to the set  $ \  \Conv_r, \ $ write $ \ \phi(\cdot) \in \Conv_r, \ $   iff the function $ \ \lambda \to \phi(|\lambda|^{1/r}) \ $ is convex. \par

\vspace{4mm}

  \ For instance, the classical subgaussian function $ \ \phi_2(\lambda) = 0.5 \lambda^2, \ \lambda \in R \ $ belongs to the set $ \ \Conv_2. \ $ \par

 \vspace{3mm}

\ It is proved  ibid that if the sequence of independent r.v. $ \ \eta_j \ $  belong to the space  $  \ B(\phi)  \ $  with $ \ \phi(\cdot) \in \Conv_2, $
then

\begin{equation} \label{ineq Pithagor}
||\sum_{j=1}^n \eta_j ||^2 B\phi \le \sum_{j=1}^n ||\eta_j||^2 B(\phi),
\end{equation}
 the Pythagoras inequality.\par

\vspace{3mm}

 \ More generally, let $ \ \phi(\cdot)  \ $ be arbitrary function from the set $ \ \Phi. \ $ Define the following  its transformation

\begin{equation} \label{hat phi}
\hat{\phi}(\lambda) \stackrel{def}{=} \sup_n \ \sup_{a \in D(n)} \phi \left( \ \lambda \sum_{j=1}^n a_j \ \right).
\end{equation}

 \ This function $ \ \hat{\phi}(\cdot) \ $  obeys a following sense. Let $ \ \{\eta_j\} \ $ be independent r.v. - s from the set $ \ B\phi \ $
and have an unit norm in this space: $ \ ||\eta_j||B\phi = 1.  \ $  Then for all the values $ \ n = 1,2,\ldots \ $

$$
||\sum_{j=1}^n  a_j \eta_j ||B\hat{\phi} \le 1, \hspace{3mm} \{a_j\} \in D(n).
$$

\vspace{5mm}

\section{Main  result. Exact Khinchine's constant calculation.}

\vspace{5mm}

 \hspace{3mm} {\bf Theorem 3.1.} Suppose that the source random variable $ \ \xi \ $ belongs to some $ \ B(\phi), \ \exists \phi \in \Phi \ $ space
 and assume besides  that the function $ \ \phi(\cdot) \ $ belongs to the class $ \ \Conv_2. \ $  Then

\vspace{4mm}

\begin{equation} \label{B result}
B[L(\xi)]\{B(\phi)\}  = ||\xi||B(\phi).
\end{equation}

\vspace{4mm}

 \ {\bf Proof } is simple.  Let $ \ \xi \in B(\phi), \ \phi \in \Conv_2. \ $ Let also $ \ \sum_{j=1}^n a^2_j = 1. \ $
We have the following upper estimate for the arbitrary integer positive value $ \ n \ $

$$
||S_n||B(\phi) \le ||\xi|| B(\phi) \cdot \sqrt{ \ \sum_{j=1}^n a^2_j  \ } = ||\xi|| B(\phi).
$$

 \ On the other hand, we have the following lower estimate, choosing the value $ \ n = 1 \ $

$$
||S_1||B(\phi) =  ||\xi|| B(\phi),
$$

 \hspace{3mm} This completes the proof of (\ref{B result}). \par

\vspace{4mm}

 \ {\bf  Example 3.1.} Let the r.v. $ \  \theta \ $ has a Rademacher's distribution; then

$$
{\bf E} e^{\lambda \theta} = \cosh \lambda, \ \lambda \in R.
$$
 \ As long as $ \ \cosh \lambda  \le \exp(\lambda^2/2), \ $ we observe that the r.v. $ \ \theta \ $ is subgaussian
and has an unit norm in the space $  \  B\phi_2. \ $ Since the function $ \ \phi_2 \ $ belongs to the set $ \ \Conv_2, \ $
one can apply the proposition of theorem 3.1. \par

\vspace{3mm}

 \ Further, the subgaussian  $ \ B\phi_2 \ $ norm of the r.v. is equivalent to the Grand Lebesgue norm with the generating function
 $ \ \psi_2(p) = \sqrt{p}, \ p \ge 1; \ $ and we obtain the known result

$$
\sup_n ||S_n[\theta]||_p \le C \ \sqrt{p},
$$
or equally

$$
\sup_n ||S_n[\theta]||B\phi_2  < \infty.
$$

\vspace{4mm}

 \ {\bf Example 3.2.} Suppose that the source  (centered) variable $ \ \nu \ $ belongs to the space $ \ B \phi_m,\  m \ge 1. \ $
 Denote $ \ m' := \min(m,2). \ $  We conclude

 $$
 \sup_n ||S_n[\nu]||B\phi_{m'} < \infty
 $$
or equally

$$
\sup_n ||S_n[\nu]||_p \le C_1(m) \ p^{1/m'}, \ p \ge 1.
$$

\vspace{4mm}

 \  Let us consider a more general case of arbitrary moment generating function $ \ \phi(\cdot) \in \Phi. \ $ \par

 \ {\bf Theorem 3.2.}

\begin{equation} \label{B general result}
B[L(\xi)]\{ B(\hat{\phi})\}  \le ||\xi||B{\phi}.
\end{equation}

\vspace{4mm}

 \ {\bf Proof } follows immediately from the direct definition of $ \ \hat{\phi}. \ $ Indeed,
    let $ \ \xi \in B(\phi), \ \phi \in \Phi. \ $ Let also as above  $ \ \sum_{j=1}^n a^2_j = 1. \ $
We have the following upper estimate for the arbitrary integer positive value $ \ n \ $

$$
||S_n||B(\hat{\phi}) \le ||\xi|| B(\phi) \cdot \sqrt{ \ \sum_{j=1}^n a^2_j  \ } = ||\xi|| B(\phi).
$$

 \vspace{5mm}

\section{Case of non - identical distributed variables.}

\vspace{5mm}

 \hspace{3mm}  It is no hard to generalize obtained  before results into the case of the non - identical distributed variables. \par

 \vspace{3mm}

 \ Let now the r.v. - s $ \  \{ \xi_k \}, \ k = 1,2,\ldots \ $ be a sequence of independent centered
 but not necessary identical distributed random variables. Suppose that the each r.v.   $ \ \xi_k \ $
 belongs to some  space $ \ B\phi_k, \ \phi_k \in \Phi: \ $

$$
{\bf E} \exp(\lambda \ \xi_k) \le \exp(\phi_k(\lambda)), \ |\lambda| \le \lambda_0, \ \exists \lambda_0 > 0.
$$

 \ Of course, one can choose all the functions $ \ \{  \  \phi_k(\cdot) \ \}, \ k = 1,2,\ldots \ $ as a natural ones
for the r.v. - s $ \ \xi_k. \ $ \par

\vspace{3mm}

 \ As before,  $ \  S_n := \sum_{k=1}^n a_k \ \xi_k.  \ $ Define the following function

\begin{equation} \label{kappa fun}
\kappa(\lambda) \stackrel{def}{=}  \sup_n \sup_{a \in D(n)} \sum_{k=1}^n \phi_k (a_k \ \lambda).
\end{equation}

 \vspace{4mm}

 \ {\bf Theorem 4.1.}   Assume that $ \  \kappa(\cdot) \in \Phi.   \  $ Then

\begin{equation} \label{not ident}
B[L\{\xi_k\}]{G\kappa} \le 1.
\end{equation}

 \vspace{4mm}

  \hspace{3mm} {\bf Proof.}  Indeed,  we have for arbitrary value $ \ a \in D(n) \ $

$$
{\bf E} \exp(\lambda S_n)  = \prod_{k=1}^n  {\bf E} \exp(\lambda a_k \ \xi_k) \le
$$

$$
\prod_{k=1}^n \exp(\phi_k(a_k \lambda))  \le \exp(\kappa(\lambda)),
$$
or equally

$$
||S_n||B(\kappa) \le 1
$$
uniformly in $ \ a, n. \ $ \par

 \vspace{3mm}

  \ {\bf Remark 4.1.} The proposition of theorem 3.1 is a particular case of considered here, indeed,
when the functions $ \ \phi_k(\cdot) \ $ are equal: $ \ \phi_k(\lambda) = \phi(\lambda). \ $ \par
 \ Therefore, it is also essentially non - improvable. \par

\vspace{5mm}

 \section{Grand Lebesgue Space approach.}

\vspace{5mm}

 \hspace{3mm}  We suppose now only that the centered  r.v. $ \ \xi \  $ belongs to certain Grand Lebesgue Space (GLS)
$ \ G\psi, \ \psi \in \Psi. \ $  As above, $ \  \{ \xi_j \} \ $ are independent copies $ \ \xi. \ $  \par

\vspace{3mm}

 \ We will apply the famous Rosenthal's  \cite{Rosenthal} inequality  for the variable $ \ C(p), \ p \ge 2, \ $ where

\begin{equation} \label{Cp}
C(p) \stackrel{def}{=}  \sup_n \sup_{ \ \{\eta_j \} }  \ \frac{||\sum_{j=1}^n \eta_j||_p}{\max \left(||\sum_{j=1}^n \eta_j||_2, \ \left( \sum_{j=1}^n ||\eta_j||_p^p  \ \right)^{1/p} \right)} < \infty.
\end{equation}

\vspace{3mm}

 \ Here $ \ \{\eta_j\}, \ j = 1,2,\ldots \ $ are independent copies of the centered variable $ \ \eta = \eta_1 \ $  such that $ \ \eta \in L(p,\Omega). \ $

\vspace{3mm}

 \ Indeed, there are huge numbers of works devoted to evaluate of these "constants", see e.g.  \cite{Dharmadhikari Jogdeo}, \cite{Ibragimov1}, \cite{Ibragimov2},
\cite{Naimark Ostrovsky}, \cite{Rosenthal} etc. In the article \cite{Naimark Ostrovsky} was obtained the  ultimate optimal {\it order} of the value of $ \ C(p): \ $

\begin{equation} \label{fun gp}
C(p) \le C_R \frac{p}{e \ \ln p}, \ p \ge 2,
\end{equation}
where

\begin{equation} \label{optimal Cp}
C_R \approx 1.776379 < 1.77638.
\end{equation}

 \ Note that for the symmetrical distributed r.v. $ \  \{ \ \eta_j \ \} $ this constant is equal (approximately)  to $ \ 1.53572.\ $ \par

\vspace{3mm}

 \ {\bf Remark 5.1.} Both the last estimates  are attainable: in the first case when the  r.v. $ \ \eta  \ $  has a form $ \ \eta = \rho - 1, \ $
 the r.v. $ \ \rho \ $ has a standard  Poisson distribution with parameter 1, \ in the second case $ \ \eta \ $ has a form $ \ \eta = \rho_1 - \rho_2, \ $
 where $ \ \rho_1, \rho_2 \ $ are independent Poisson distributed with parameter  0.5; \ see  \cite{Ibragimov1}, \cite{Ibragimov2}, \cite{Naimark Ostrovsky}. \par

 \vspace{3mm}

  \ Assume now that the source  centered r.v. $ \ \xi \ $ belongs to some  GLS $ \ G\psi, \ \psi \in \Psi. \ $ Introduce the auxiliary such a function

 \begin{equation} \label{auxil}
 \psi_R(p) := C_R \ \frac{p}{e \ \ln p} \ \psi(p), \ p \ge 2.
 \end{equation}

\vspace{4mm}

 \ {\bf Theorem 5.1.}

\vspace{3mm}

\begin{equation} \label{Gpsi est}
B[L(\xi)](G\psi_R) \le ||\xi||G\psi.
\end{equation}

\vspace{4mm}

 \ {\bf Proof.} Let $ \ \xi \in G\psi; \ $ one can suppose without loss of generality $ \ ||\xi||G\psi = 1. \ $ Therefore, for all the  acceptable values
  $ \ p \  \Rightarrow \ ||\xi||_p \le \psi(p). \ $ \par

  \ We intent to apply the estimate  (\ref{Cp}) with clarification  (\ref{optimal Cp}), substituting the r.v.$ \ a_j \ \xi_j \ $
  instead $ \ \eta_j. \ $ We have denoting $ \ \sigma^2 = \Var(\xi), \ $ taking into account the equality $ \ a \in D(n): \ $

 $$
 ||\sum_{j=1}^n a_j \xi_j||^2_2 = \sigma^2;
 $$

$$
\sum_{j=1}^n |a_j|^p \ ||\xi_j||^p_p \le \psi^p(p) \ \sum_{j=1}^n |a_j|^p = \psi^p(p) \ ||a||_p^p \le \psi^p(p) \ ||a||_2^p,
$$
as long as for $ \ p \ge 2 \ \Rightarrow ||a||_p \le ||a||_2 = 1. \ $ Note yet $ \  \sigma \le \psi(p), \ $ as long as $ \ p \ge 2. \ $ \par

 \ We obtain by virtue of (\ref{Cp})

$$
||\sum_{j=1}^n \ a_j \xi_j||_p \le C(p) \max \left[ \  \sigma \ [\sum_{j=1}^n a_j^2]^{1/2}, \ ||a||_p \ \psi(p) \ \right] \le
$$

$$
C_R \ \frac{p}{e \ \ln p} \ \psi(p) = \psi_R(p),
$$
Q.E.D. \par

\vspace{5mm}

\section{Khintchine's inequality in the space of continuous functions.}

\vspace{5mm}

 \hspace{3mm} Let $ \  Z = \{z\}  \ $ be arbitrary set; the semi - distance function  $ \  \rho = \rho(z_1,z_2), \ z_{1,2} \in Z \ $ on this set will be
 introduced below. Recall that the semi - distance function is non - negative symmetrical function vanishing in the diagonal $ \ \rho(z,z) = 0, \ $
 satisfying the triangle inequality but  in general case the relation $ \ \rho(z_1,z_2) = 0 \ $ does not imply that  $ \ z_2 = z_1. \ $ \par
 The (Banach) space of all the continuous numerical valued functions will be denoted as ordinary $ \ C(Z) = C(Z,\rho); $ it is equipped with the uniform
 norm

 $$
 ||f||C(Z) := \sup_{z \in Z} |f(z)|.
 $$

 \vspace{3mm}

 \hspace{3mm} Let  $ \  \eta = \eta(z),  \ z \in Z  \ $ be separable centered: $ \ {\bf E} \eta(z) = 0 \ $ numerical valued random field (r.f.).
 Let also $ \ \eta_i = \eta_i(z), \ i = 1,2,\ldots \ $ be independent copies of $ \ \eta(z). \ $
  \ Let us impose the following  {\it subgaussian } condition on the r.f. $ \ \eta = \eta(z): \ $

\begin{equation} \label{unif subgaussian}
\sigma :=   \sup_{z \in Z} ||\eta(z)||B(\phi_2)  < \infty.
\end{equation}
 \ Introduce therefore the following  bounded {\it semi - distance } on  the set $ \ Z: \ $

 $$
\rho(z_1,z_2) \stackrel{def}{=} ||\eta(z_1) - \eta(z_2)||B(\phi_2).
 $$

 \ Denote also by $ \ H(\epsilon) = H(Z,\rho,\epsilon), \ 0 < \epsilon \le C_5, \ $ the metric entropy of the whole set (space) $ \ Z \ $
relative the metric $ \ \rho, \ $  i.e. the  (natural) logarithm of the minimal amount of the  closed ball in the distance $ \ \rho, \ $
which cover this set $ \ Z. \ $

\vspace{4mm}

  \ Put as above for arbitrary tuple $ \ a = \vec{a}(n) \in D(n) \ $

 $$
 Y_{a(n)}(z) \stackrel{def}{=}  \sum_{i=1}^n  a_i \eta_i(z), \ a(n) = \vec{a(n)} \in D(n),
 $$
and define the r.v.

$$
\beta \stackrel{def}{=} \sup_{a(n) \in D(n)} ||\sup_{z \in Z}   \ Y_{a(n)} (z)|| B(\phi_2).
$$

\vspace{5mm}

 \ {\bf Theorem 6.1.} Suppose that the following entropic integral convergent:

\begin{equation} \label{entropic integral}
\int_0^1 H^{1/2}(Z,\rho, \epsilon) \ d \epsilon < \infty \  -
\end{equation}
\ the famous Dudley condition, \cite{Dudley}, \cite{Dudley CLT}. \ Then all the random fields $ \ Y_{a(n)}(z) \ $
are $ \ \rho \ - \ $ continuous with probability one:

$$
{\bf P} \left( \ Y_{a(n)}(\cdot) \in C(Z,\rho) \ \right) = 1
$$
 and the r.v. $ \ \beta \ $  is subgaussian: $ \ \beta \in B(\phi_2) \ $
or equally:

$$
\exists C_0 = C_0(H, \sigma) = \const < \infty \ \Rightarrow ||\beta||_p \le C_0  \ \sqrt{p} < \infty.
$$

\vspace{4mm}

 \ {\bf Proof.} Let us consider the  subgaussian  random field $ \  Y_{a(n)}(z). \ $  We have taking into account the equality
  $ \   \sum_i a^2(i) = 1 \ $ and properties of the subgaussian norm

$$
\sup_{z \in Z} \ \sup_n \max_{a(n) \in D(n)} ||Y_{a(n)}(z)||B(\phi_2) \le \sigma
$$
and

$$
 \sup_n \max_{a(n) \in D(n)} ||Y_{a(n)}(z_1)   -  Y_{a(n)}(z_2)||B(\phi_2) \le \rho(z_1, z_2).
$$

 \ Both the propositions of theorem 6.1  follows immediately from Theorem 3.17.1
of monograph   \cite{Ostrovsky1}, chapter 3, section 17; see also \cite{Buldygin 2}, \cite{Buldygin 3}. \par

\vspace{5mm}

 \section{Concluding remarks.}

\vspace{5mm}

  \ It is interest in our opinion to generalize the obtained results on the sequence of
  the centered (stationary or not)  random variables satisfying one or another mixing condition. \par

\vspace{6mm}

\vspace{0.5cm} \emph{Acknowledgement.} {\footnotesize The first
author has been partially supported by the Gruppo Nazionale per
l'Analisi Matematica, la Probabilit\`a e le loro Applicazioni
(GNAMPA) of the Istituto Nazionale di Alta Matematica (INdAM) and by
Universit\`a degli Studi di Napoli Parthenope through the project
\lq\lq sostegno alla Ricerca individuale\rq\rq .\par

\end{document}